\documentclass[12pt]{article}

\usepackage[letterpaper,margin=1in]{geometry}
\usepackage[T1]{fontenc}
\usepackage{amsmath,amssymb,amsthm,mathtools}
\usepackage{newpxtext,newpxmath}
\usepackage{microtype}
\usepackage{graphicx}
\usepackage{xcolor}
\usepackage{tikz}
\usetikzlibrary{arrows.meta,calc,positioning}
\usepackage[hidelinks]{hyperref}
\usepackage{enumitem}
\usepackage{caption}
\usepackage{titlesec}

\titleformat{\section}{\Large\bfseries}{\thesection.}{0.55em}{}
\titleformat{\subsection}{\large\bfseries}{\thesubsection.}{0.5em}{}
\newtheorem{theorem}{Theorem}
\newtheorem{proposition}[theorem]{Proposition}
\newtheorem{lemma}[theorem]{Lemma}

\theoremstyle{definition}
\newtheorem*{problem}{Problem}

\theoremstyle{remark}

\newcommand{\mG}{\mathbf G}
\newcommand{\Z}{\mathbb Z}

\newcommand{\OO}{\mathcal O}
\newcommand{\Ap}{\mathcal A}

\newcommand{\symb}[2]{\left(\frac{#1}{#2}\right)}
\newcommand{\SL}{\operatorname{SL}}

\providecommand{\doi}[1]{\href{https://doi.org/#1}{doi:\,#1}}

\title{Squares, Three Fleas, Sporadic Sets, and Squares}
\author{Giedrius Alkauskas}
\date{}

\begin{document}
\maketitle

\begin{abstract}
Three fleas are on the plane. At each move, two fleas
jump to the other two vertices of a square erected on the segment joining
them. Starting from $\mG=\{(0,0),(2,1),(3,2)\}$, we ask which areas can occur for the triangle they span. Every positive half-integer occurs, and the attainable integer areas meet every residue class modulo every modulus. Yet no perfect square ever occurs. We explain this unexpected obstruction using quadratic reciprocity, and then uncover its geometric source: the flea dynamics are closely related to integral Apollonian circle packings. An exact computation up to \(10^{12}\) finds only five additional missing nonsquares: $5,29,80,99,179$.
\end{abstract}

\section{Geometry}
How far beyond olympiad mathematics can an olympiad problem lead? Consider one example \cite{mse}.
\begin{problem}
In the plane, three fleas are given. During a jump, two fleas at points
$A$ and $B$ move to the other two vertices $C,D$ of a square $ABCD$, with
$AB$ as a side. Suppose that no possible first jump decreases the area of
the triangle spanned by the fleas. Prove that no finite sequence of jumps
can decrease that area.
\end{problem}

\begin{figure}[h]
  \centering
  \includegraphics[width=0.68\textwidth]{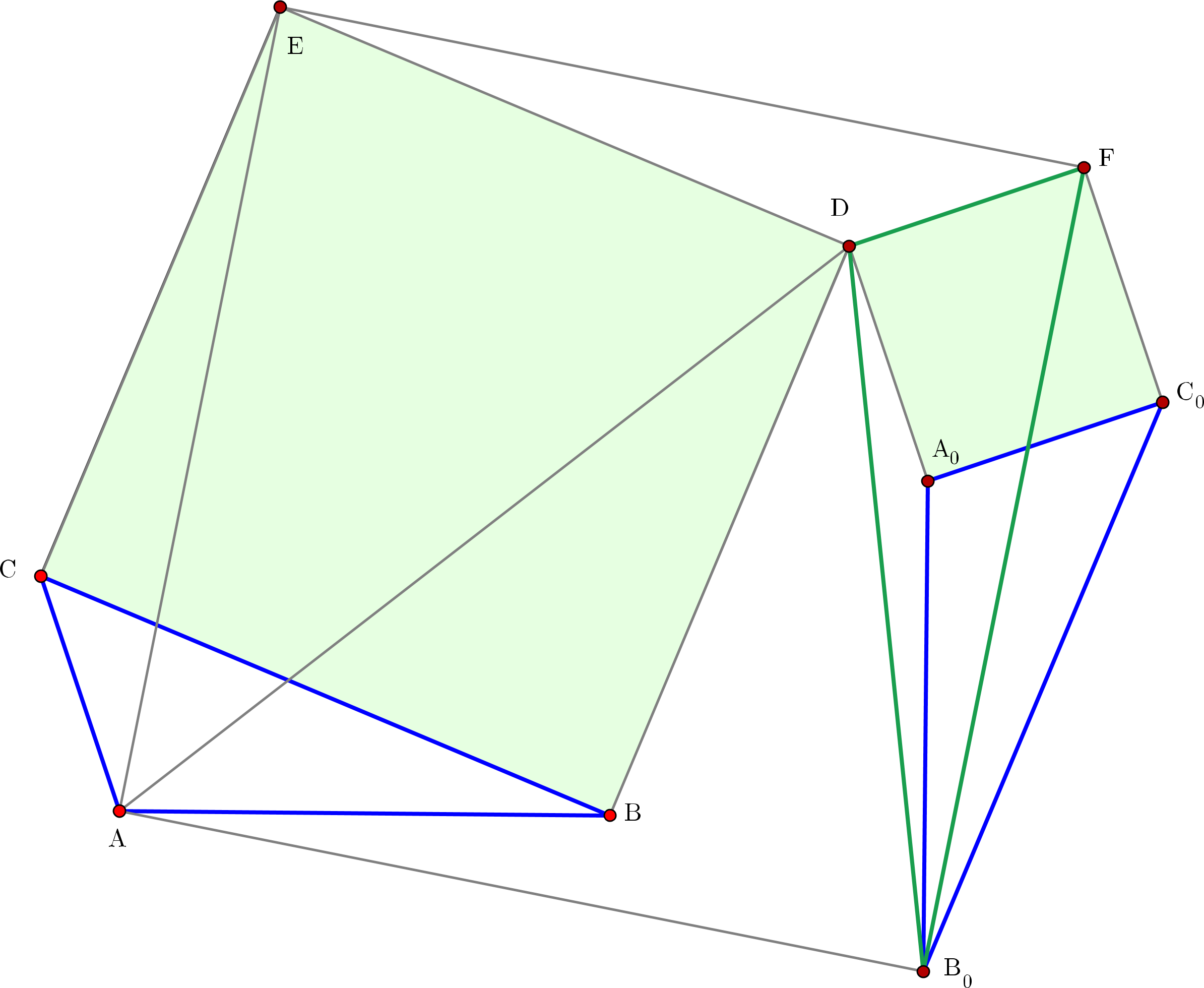}
  \caption{Two jumps of opposite signs can be replaced by one jump.}
  \label{cong}
\end{figure}

Here is the central idea of the solution. Once the three fleas are cyclically ordered, the two
squares on a chosen side define a \emph{positive} and a \emph{negative}
jump. Fix an initial triangle $\mathcal T$ and a target reachable triangle
$\mathcal V$, and choose a shortest sequence of jumps from $\mathcal T$ to
a triangle congruent to $\mathcal V$. Such a sequence cannot change sign.
Indeed, consecutive jumps of opposite signs cancel if they use the same side; otherwise the geometric identity in  Figure~\ref{cong} replaces them by one jump, up to relabelling. A mixed shortest sequence would therefore not be shortest.

Orient the initial triangle positively. Positive jumps increase its oriented
area. If the first negative jump preserved the orientation, it would reduce
the unsigned area, contrary to the hypothesis; hence it must reverse the
orientation. Every subsequent negative jump then makes the signed area more
negative and so increases its absolute value. Since a shortest sequence has
only one sign, the problem follows.

\begin{figure}[h]
  \centering
  \includegraphics[width=0.62\textwidth]{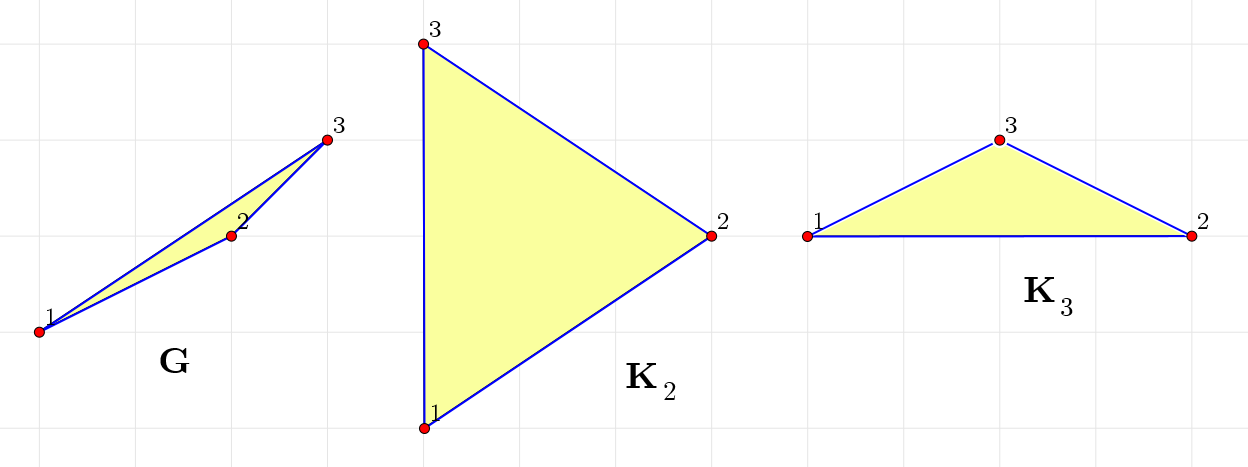}
  \caption{The three configurations obtained from $\mG$ by a negative first
  jump, up to congruence and relabelling.}
  \label{trik}
\end{figure}

Now consider the three fleas at $\{(0,0),(2,1),(3,2)\}$, which we call the
configuration $\mG$ (Figure~\ref{trik}, left). Its initial area is $1/2$,
and its three squared side lengths are $2,5,13$. A direct check of the three
possible negative jumps shows that the resulting triangle is, up to
congruence and relabelling, either $\mG$ itself or one of two isosceles
triangles, denoted by $\mathbf K_2$ and $\mathbf K_3$. Every positive
jump increases the oriented area. Thus the problem hypothesis holds, and
in particular the fleas can never become collinear.

Chronologically, this special problem for $\mG$ occurred first. However, if one tried to rule out area $0$ for $\mG$ by searching for a simple invariant of all attainable areas, an attempt would be futile at olympiad level. How can this be, if a logically stronger statement is eligible for a contest? Because a direct attack on $\mG$ might push a student immediately into investigation of arithmetic, and this is a trap.

To understand not merely whether the area decreases, but which areas occur at all, it is useful to encode flea jumps algebraically. Take the following three $\SL_4(\Z)$ matrices on faith; Section 2 will explain where they come from:
\begin{equation}\label{eq:UVW}
U=\begin{pmatrix}
1&1&1&2\\ 0&1&0&0\\ 0&0&1&0\\ 0&1&1&1
\end{pmatrix},\quad
V=\begin{pmatrix}
1&0&0&0\\ 1&1&1&2\\ 0&0&1&0\\ 1&0&1&1
\end{pmatrix},\quad
W=\begin{pmatrix}
1&0&0&0\\ 0&1&0&0\\ 1&1&1&2\\ 1&1&0&1
\end{pmatrix}.
\end{equation}
Write $\Gamma=\langle U,V,W\rangle$, and let $\gamma^+$ and $\gamma^-$ be
the semigroups generated by $U,V,W$ and by $U^{-1},V^{-1},W^{-1}$,
respectively. Surprisingly, the geometric argument illustrated in Figure \ref{cong} has the following algebraic counterpart.

\begin{theorem}\label{thm:normal-form}
There is a six-element subgroup $\Omega<\Gamma$, isomorphic to $S_3$, such
that every $g\in\Gamma$ has one of the forms $g=P\omega$, with
$P\in\gamma^+$ or $P\in\gamma^-$ and $\omega\in\Omega$; the word $P$ may
be empty. Both $\gamma^+$ and $\gamma^-$ are free semigroups on their
displayed generators.
\end{theorem}

The group $\Gamma$ itself -- to be called \emph{flea group} -- is certainly not free. For example, $(UV^{-1})^3=I$. As for the attainable areas, the actual picture is more involved.
\begin{theorem}\label{thm2}
For the flea orbit starting from $\mG$:
\begin{enumerate}[label=\textup{(\roman*)}]
\item every positive half-integer $k+\tfrac12$, $k\ge0$, occurs as an area;
      the attainable integer areas meet every residue class modulo every
      positive integer;
\item no perfect square, including $0$, occurs as an area;
\item among the positive integer nonsquares not exceeding $10^{12}$, the only
      missing values are $5,29,80,99,179$.
\end{enumerate}
\end{theorem}
The set of reachable areas is plainly enumerable. Moreover, it is decidable: Theorem~\ref{thm:normal-form}, together with the
strict growth of area under positive jumps, gives a finite exact procedure
to determine whether any prescribed number occurs. Concerning the sporadic values, it is safe to bet that the given list is complete. Yet, this is still unproved.

\section{Algebra}
We now explain the matrices in \eqref{eq:UVW}. Translate one flea to the
origin and write the other two position vectors as $u,v\in\Z^2$. Put
$s=\det(u,v)$, so the unsigned area is $|s|/2$, and define
\begin{equation*}
a=u\cdot v,\quad b=(u-v)\cdot u,\quad
c=(v-u)\cdot v.
\end{equation*}
Then $a+b=|u|^2$, $a+c=|v|^2$, $b+c=|u-v|^2$. Further, $\det(u,v)^2+(u\cdot v)^2=|u|^2\cdot |v|^2$. In expanded form, this is just a classical Brahmagupta–Fibonacci two-square identity 
\[
(x_1^2+y_1^2)(x_2^2+y_2^2)=(x_1x_2+y_1y_2)^2+(x_1y_2-x_2y_1)^2.
\]
Thus,
\begin{equation}\label{eq:s}
s^2=ab+ac+bc.
\end{equation}
Consider the jump on the side joining the endpoints of $u$ and $v$. The
quantities $b,c$ remain unchanged, while the oriented doubled area becomes
$s'=s\pm|u-v|^2=s\pm(b+c)$. By \eqref{eq:s},
$a'=((s')^2-bc)/(b+c)=a+(b+c)\pm2s$. Thus the two jumps on this side are
\begin{equation}\label{eq:jump}
(a,b,c;s)\longmapsto
\bigl(a+b+c\pm2s,\ b,\ c;\ s\pm(b+c)\bigr),
\end{equation}
and the other four possibilities are obtained by permuting $a,b,c$. The
plus sign in \eqref{eq:jump} is exactly the action of $U$ in
\eqref{eq:UVW} on the column vector $(a,b,c,s)^T$; similarly, $V,W$ are the
other two positive jumps. More generally,
\begin{equation}\label{power}
	U^t(a,b,c;s)=
	\bigl(a+(b+c)t^2+2st,\ b,\ c;\ s+t(b+c)\bigr),\qquad t\in\Z,
\end{equation}
with analogous formulas for $V^t,W^t$. Some relations among these matrices have a simple geometric meaning. Put
\[
T_1=\begin{pmatrix}
1&0&0&0\\0&0&1&0\\0&1&0&0\\0&0&0&-1
\end{pmatrix},\quad
T_2=\begin{pmatrix}
0&0&1&0\\0&1&0&0\\1&0&0&0\\0&0&0&-1
\end{pmatrix},\quad
T_3=\begin{pmatrix}
0&1&0&0\\1&0&0&0\\0&0&1&0\\0&0&0&-1
\end{pmatrix}.
\]
Each $T_i$ transposes two of $a,b,c$ and reverses orientation. These belong to $\Gamma$ since, for example, $T_{3}=UV^{-1}U$. Furthermore,
\begin{equation}\label{eq:T1-conjugation}
T_1UT_1=U^{-1},\qquad T_1VT_1=W^{-1},\qquad
T_1WT_1=V^{-1},
\end{equation}
with cyclic analogues for $T_2,T_3$. Finally,
$C=T_3T_2=U^{-1}VW^{-1}U$ cycles $a,b,c$ and has order $3$. Hence
$\Omega=\{I,T_1,T_2,T_3,C,C^2\}\cong S_3$.

We now prove Theorem~\ref{thm:normal-form}. For distinct generators, the six
mixed products reduce as follows:
\begin{equation}\label{eq:mixed-reductions}
\begin{aligned}
UV^{-1}&=VT_3,&VU^{-1}&=UT_3,&
UW^{-1}&=WT_2,\\
WU^{-1}&=UT_2,& VW^{-1}&=WT_1,& WV^{-1}&=VT_1.
\end{aligned}
\end{equation}
The corresponding identities with the signs reversed follow by inversion
and conjugation. Since $\Omega$ normalizes
$\{U^{\pm1},V^{\pm1},W^{\pm1}\}$, a factor from $\Omega$ may always be
moved to the right (geometrically, this is simply the phrase “up to relabelling the vertices” from Section 1).

Among all expressions $g=P_1^{\varepsilon_1}\cdots
P_r^{\varepsilon_r}\omega$, where $P_i\in\{U,V,W\}$,
$\varepsilon_i\in\{\pm1\}$, and $\omega\in\Omega$, choose one with
minimal $r$. If the signs are mixed, there is an adjacent pair with opposite
signs. If the underlying generators agree, they cancel; if they differ,
\eqref{eq:mixed-reductions} replaces the two letters by one signed generator
and an element of $\Omega$, which is moved to the right. In either case
$r$ decreases, a contradiction. Thus all signs agree, proving the one-sign
normal form.

Let us see the algorithm at work. Take $P=VW^{-1}U^{-2}W^3V^2$. Since
$VW^{-1}=WT_1$, equations \eqref{eq:T1-conjugation} give
\[
P=(VW^{-1}T_1)(T_1U^{-2}T_1)(T_1W^3T_1)(T_1V^2T_1)T_1
  =WU^2V^{-3}W^{-2}T_1.
\]
Now $UV^{-1}=VT_3$, and the same procedure gives
$P=WUVU^2W^2T_3T_1=WUVU^2W^2C^2$. Thus a rather mixed word has become a
positive word followed by a harmless relabelling.

To prove freeness of a positive semigroup, we use the standard ping-pong argument. Consider the open cone $\mathcal C=\{(a,b,c;s)\in\mathbb R^4:a,b,c,s>0\}$. All three
generators map $\mathcal C$ into itself, and
\[
U\mathcal C\subset\{a>\max(b,c)\},\qquad
V\mathcal C\subset\{b>\max(a,c)\},\qquad
W\mathcal C\subset\{c>\max(a,b)\}.
\]
These regions are disjoint. First apply a word to $(1,1,1;1)$: the empty
word leaves this point outside all three regions, whereas every nonempty
word lands in the region determined by its leftmost letter. Thus no
nonempty positive word is the identity. If two nonempty positive words
define the same matrix, their leftmost letters must agree; cancel them and
continue by induction. Hence the words are identical. Taking inverses proves
the same for $\gamma^-$. This completes the proof of
Theorem~\ref{thm:normal-form}.

The algebra now pays for itself: it turns a two-sided group orbit into three
monotone trees. In flea coordinates, $\mG=(8,-3,5;1)$. As noted, each negative first jump changes the orientation. Reversing it back and relabelling the fleas gives
\begin{equation*}
	K_1=(8,5,-3;1),\qquad K_2=(5,8,8;12),\qquad
	K_3=(-3,8,8;4).
\end{equation*}
For a positive state $P=(a,b,c;s)$ (so, $s>0$), let
$\mathscr S(P)$ denote all doubled areas $s$ attainable from $P$ by applying only positive jumps. Accordingly, we have a decomposition
\begin{equation*}
	\{|s(P)|:P\in\Gamma\mG\}
	=\mathscr S(\mG)\cup\mathscr S(K_2)\cup\mathscr S(K_3).
\end{equation*}
As for the half-integers, they need no search: formula \eqref{power} gives
$s(U^t\mG)=1+2t$ for $t\ge0$. Hence every $t+\tfrac12$ occurs. Further, 
every positive jump increases $s$ by a positive squared side length. Thus,
to decide whether an integral area occurs, only finite search is needed. It
certifies that up to $10^{12}$ only squares and $5,29,80,99,179$ are missed. The absence of these five values is, certainly, quaint. But why are squares missed?

\section{Arithmetic}
First, the missing squares are not the result of an ordinary congruence sieve.
\begin{proposition}\label{prop:all-residues}
For every $m\ge1$ and every residue $r\pmod m$, some reachable integer area
is congruent to $r\pmod m$.
\end{proposition}
\begin{proof}
Apply $V^t$ to $\mG$ and then $U^{2k}$, with $t$ odd.
Formulas \eqref{power} give the integer area
\begin{equation}
\frac{1+13t}{2}+k(13t^2+2t+2).
\label{ars}
\end{equation}
For every odd prime $p\mid m$, the quadratic $13t^2+2t+2$ is not identically
zero modulo $p$. So, for a certain $t$, it is coprime to $p$. Further, $13t^2+2t+2$ is always odd for odd $t$. The Chinese Remainder Theorem (CRT) therefore supplies a
positive odd $t$ for which $13t^2+2t+2$ is coprime to $m$. With this $t$
fixed, varying $k\ge0$ reaches every residue modulo $m$.
\end{proof}

Thus no finite collection of congruence conditions can explain the missing
squares. What, then, can? To predict the answer, consider an example. Equation \eqref{ars} (for $t=2r+1$) gives $7+13r+k(52r^2+56r+17)$ which is never a square for $r,k\in\mathbb{N}$. For $r=2$, this reduces to the fact that $\symb{33}{337}=-1$, in Legendre symbol notation. Yet, for $r=1$, the set specialises to $20+125k$. Again no square occurs, now because $v_{5}(20+125k)=1$.

 Experiments quickly produce different local arguments of this kind. It appears that, provided elementary arithmetic is handled with care, the uniform explanation comes from quadratic reciprocity. The proof has two halves. First, a sign equal to $-1$ of a suitably defined quadratic character travels with every odd coordinate in the orbit of $\mG
$. Second, a square area would force that same sign to be $+1$. We use $\symb{u}{p}$ for the Legendre symbol when $p$ is an odd prime, and its Jacobi or Kronecker extension when the denominator is composite or even. 

Set $d=a+b+c-2s$. Equation \eqref{eq:s} becomes
\begin{equation}\label{eq:Descartes}
(a+b+c+d)^2=2(a^2+b^2+c^2+d^2).
\end{equation}
A connoisseur will immediately recognize Descartes equation for curvatures of $4$ mutually touching circles, though there are no circles in flea geometry! The character argument below is a self-contained specialization of the quadratic-reciprocity mechanism of Haag, Kertzer, Rickards, and Stange \cite{HKRS}. Lemma \ref{lem:ar} supplies the additional step needed to treat flea areas rather than circle curvatures.

For a quadruple $x=(x_1,x_2,x_3,x_4)$ satisfying
\eqref{eq:Descartes}, let $\sigma_i$ replace $x_i$ by
$2\sum_{j\ne i}x_j-x_i$ and fix the other three entries (this is standard replacement of one root of the quadratic equation with another one). In the variables
$(a,b,c,d)$, a positive $U$-jump is $\sigma_4$ followed by the transposition
of the first and fourth coordinates. Hence every flea state from $\mG$ lies
in the enlarged orbit $\OO$ of $Q_0=(-3,5,8,8)$ under the four reflections
$\sigma_i$ and coordinate permutations.

\begin{lemma}\label{lem1}
Every $Q=(x_1,x_2,x_3,x_4)\in\OO$ is primitive. Modulo $8$, two entries are
$5$ and the other two lie in $\{0,4\}$. Moreover, every pair sum
$x_i+x_j$, $i\neq j$, is positive. For every odd prime $p$, at most two of the
coordinates $x_i$ are divisible by $p$.
\end{lemma}
\begin{proof}
Primitivity is preserved because each $\sigma_i$ is an integral involution.
For the residue pattern, a $5$ is replaced by a $5$, while an even residue
$e\in\{0,4\}$ is replaced by $4-e$. Finally, every partition
$\{i,j,k,\ell\}=\{1,2,3,4\}$ satisfies
\begin{equation*}
2(x_i+x_j)(x_k+x_\ell)=(x_i-x_j)^2+(x_k-x_\ell)^2.
\end{equation*}
All pair sums are positive for $Q_0$. Under a reflection the complementary
pair sum is unchanged, the above identity prevents a new pair sum from
being negative. Equality to $0$ and primitivity would force, up to permutation, the
quadruple $(0,0,1,1)$, contrary to the residue pattern. The last statement for odd $p$ follows immediately from \eqref{eq:Descartes}.
\end{proof}
Thus every $Q\in\OO$ has two odd entries, both $5\pmod8$, and two entries
divisible by $4$.  We shall now attach a sign to either odd coordinate. That sign can be read locally modulo every prime dividing that coordinate, yet survives every Descartes reflection. Mark an odd entry $n$. For every odd prime $p\mid n$,
choose any other coordinate $m\not\equiv0\pmod p$ (at least one such exists) and put
$\varepsilon_p(Q,n)=\symb{m}{p}$. This is well defined: modulo $p$, the
other three entries $x,y,z$, according to \eqref{eq:Descartes}, satisfy $(x-y-z)^2\equiv4yz$ and its cyclic
analogues, so all their nonzero members have the same Legendre symbol. Let
\begin{equation}\label{eq:character}
\chi_Q(n)=\prod_{\substack{p\mid n\\v_p(n)\text{ odd}}}
\varepsilon_p(Q,n).
\end{equation}
If another coordinate $m$ is coprime to $n$, then
$\chi_Q(n)=\symb{m}{|n|}$. If $n$ is held fixed while another coordinate is
reflected, the character is unchanged: for each $p\mid n$, one of the two
unchanged coordinates is (Lemma \ref{lem1}) nonzero modulo $p$ and may be used before
and after the reflection. We will now show that \emph{both} characters remain unchanged.

\begin{lemma}\label{lem:two-odd}
The two odd entries of every $Q\in\OO$ have the same constant character on $\OO$:
\begin{equation*}
\chi(Q)=-1\qquad(Q\in\OO).
\end{equation*}
\end{lemma}

\begin{proof}
Write $Q=(a,b,c,d)$, with $a,b$ odd and $c,d$ even. Keeping $a,b$ fixed and
alternately reflecting the last two coordinates produces
$f(t)=(a+b)t^2-(a+b+c-d)t+c$, and, up to interchanging the last two entries,
$(a,b,f(t),f(t+1))\in\OO$ for every $t\in\Z$. For every prime $p\mid ab$, the polynomial $f$ is not identically zero modulo $p$. Indeed, if the first two coefficients vanish modulo $p$, then $p\nmid c$. In other cases, one can also choose \(t\) with \(p\nmid f(t)\). CRT therefore gives $t$ such that $e=f(t)$,
$\gcd(e,ab)=1$. Since $a+b>0$, choosing a sufficiently large representative, we may arrange $e>0$. This yields the orbit point $(a,b,e,h)$. Also, $4\mid e,h$ by Lemma \ref{lem1}. Now,
\[
\chi_Q(a)=\symb{e}{|a|},\qquad \chi_Q(b)=\symb{e}{|b|}.
\]
If $u$ is odd and $u\equiv1\pmod4$ as a signed integer, while $4\mid e$ and
$\gcd(u,e)=1$, quadratic reciprocity and the supplementary law for $2$ give
$\symb{e}{|u|}=\symb{u}{e}$. Finally,  \eqref{eq:Descartes} yields
\[
\left(\frac{a+b-e-h}{2}\right)^2=ab+eh.
\]
Thus $ab$ is an invertible square modulo $e$, so
$1=\symb{ab}{e}=\chi_Q(a)\chi_Q(b)$. The two signs are therefore equal. Under any reflection at least one odd entry is unchanged, so the other remains unchanged, too. Consequently, the persisting common sign is the same as the root. Namely, $\symb{8}{5}=-1$.
\end{proof}

The next lemma is the flea-specific bridge.

\begin{lemma}\label{lem:ar}
	Suppose that a flea state $(a,b,c;s)$ has primitive associated Descartes
	quadruple $Q$ and that $|s|=2w^2$, $w\ge1$. If exactly one of $a,b,c$ is odd,
	denote it by $n$. Then $\chi_Q(n)=1$.
\end{lemma}

\begin{proof}
	Relabel so that $n=a$, and put $R=a+b$, $T=a+c$. 
	The squared side lengths \(R=a+b\) and \(T=a+c\) are positive and, since \(b,c\) are even, odd. Now,
	\[
	RT=(a+b)(a+c)=a^2+s^2=n^2+4w^4.
	\]
Choose (as allowed by CRT) an integer $z$ such that	$z\equiv w\pmod{R^2}$, $z\equiv1\pmod p$ for every prime $p\mid n$ with $p\nmid R$. We shall manufacture an odd integer \(E\), coprime to \(n\), which reproduces the local character of \(n\) at every prime dividing \(n\), but whose own prime divisors force \(n\) to be a quadratic residue. Thus, define
	\[
	E=\frac{n^2+4z^4}{R}.
	\]
This is a positive odd integer: indeed,	$n^2+4z^4=RT+4(z^4-w^4)$, and $R\mid z^4-w^4$. We first identify the local character of $n$. Let $p\mid n$. If $p\nmid R$, then $E\equiv4R^{-1}\pmod p$, and therefore $\symb{E}{p}=\symb{R}{p}$. Since $R=n+b\equiv b\pmod p$, this is precisely the local character of $n$.
	
Suppose instead that $p\mid R$. From $RT=n^2+4w^4$ we obtain $p\mid w$. Moreover $p\nmid T$: otherwise $p$ would divide $a,b,c,s$, and hence also the fourth Descartes coordinate, contrary to primitivity. Since $z\equiv w\pmod{R^2}$,
	\[
	E-T=4\,\frac{z^4-w^4}{R}
	\]
is divisible by $R$, and hence $E\equiv T\pmod p$. Thus $\symb{E}{p}=\symb{T}{p}$, again the local character of $n$, since $T=n+c\equiv c\pmod p$. In either case $p\nmid E$. Hence $\gcd(E,n)=1$, and the definition of the character gives
	\begin{equation}
	\chi_Q(n)=\symb{E}{|n|}.                     
	\label{reco}
	\end{equation}
	It remains to evaluate this Jacobi symbol. Let $q$ be any prime divisor
	of $E$. Then $q\nmid n$, and $n^2+4z^4\equiv0\pmod q$. Also $q\nmid z$, for otherwise $q\mid n$. Consequently
	\[
	i=\frac{|n|}{2z^2}\pmod q
	\]
is defined and satisfies $i^2\equiv-1\pmod q$. Thus $q\equiv1\pmod4$. Moreover, since $(1+i)^2=2i$, $|n|\equiv2iz^2=\bigl((1+i)z\bigr)^2\pmod q$, and therefore $\symb{|n|}{q}=1$. Every prime divisor of the odd integer $E$ is therefore $1\pmod4$, so $E\equiv1\pmod4$. Quadratic reciprocity and \eqref{reco} now give
	\[
	\chi_Q(n)
	=\symb{E}{|n|}
	=\symb{|n|}{E}
	=1.
	\]
\end{proof}
To finish the proof of Theorem \ref{thm2}, assume $|s|=2w^2$, $w\geq 1$. From \eqref{eq:s} and Lemma \ref{lem1}, we see that exactly one of $a,b,c$ is odd. Hence Lemmas \ref{lem:two-odd} and \ref{lem:ar} apply, and they force a contradiction.
\section{Apollonian circles}\label{relations}

Equation~\eqref{eq:Descartes} is Descartes' circle equation \cite{Bradford,levrie}. Its four
entries may be read as the signed curvatures, or \emph{bends}, of four
mutually tangent oriented circles: a negative bend denotes an enclosing
circle, and bend \(0\) denotes a line. 

For a flea triangle this interpretation is more than a formal consequence
of one quadratic identity. The associated circle configuration can also be
recovered geometrically from the triangle, for instance by using complex
coordinates; compare the complex forms of Descartes' theorem in
\cite{northshield,tupan}. We leave that construction to the interested
reader. For us, the striking point is simpler: Descartes' equation has
emerged from nothing more than the dot-product identity of Section~2, while
the Apollonian replacement of one tangent circle by another is mirrored by
an elementary jump across a square. The geometry has become elementary; the
arithmetic has not.

\begin{figure}[ht]
	\centering
	\begin{tikzpicture}[scale=.76,line cap=round,line join=round]
		\begin{scope}[xshift=-4.55cm,yshift=-.25cm,scale=.86]
			\node[font=\small] at (1.45,3.05) {the flea triangle $\mG$};
			\draw[step=1cm,gray!18,very thin] (-.35,-.35) grid (3.35,2.35);
			\coordinate (P1) at (0,0);
			\coordinate (P2) at (2,1);
			\coordinate (P3) at (3,2);
			\draw[line width=1.25pt,blue!68!black] (P1)--(P2)--(P3)--cycle;
			\foreach \P in {P1,P2,P3}
			\node[circle,fill=black,inner sep=1.65pt] at (\P) {};
			\node[font=\scriptsize,below left=2pt and 2pt of P1] {$(0,0)$};
			\node[font=\scriptsize,below right=1pt and 1pt of P2] {$(2,1)$};
			\node[font=\scriptsize,right=3pt of P3] {$(3,2)$};
			\node[font=\scriptsize,fill=white,inner sep=1pt] at (1.08,.35) {$5$};
			\node[font=\scriptsize,fill=white,inner sep=1pt] at (2.67,1.28) {$2$};
			\node[font=\scriptsize,fill=white,inner sep=1pt] at (1.38,1.18) {$13$};
		\end{scope}
		\draw[-{Stealth[length=3mm]},line width=.8pt] (-.70,.45)--(.72,.45)
		node[midway,above,font=\small] {$(a,b,c,d)$};
		\begin{scope}[xshift=4.45cm,yshift=.85cm,scale=.55]
			\node[font=\small] at (-1.25,4.05) {the Descartes configuration};
			\coordinate (Za) at (0,0);
			\coordinate (Zb) at (-1.25,-1.6666667);
			\coordinate (Zc) at (-1.25,-3);
			\coordinate (Zd) at (-2.5,0);
			\draw[line width=1.05pt,gray!72,dashed] (Zb) circle (3.333333);
			\draw[line width=1.2pt,green!45!black,fill=green!5] (Za) circle (1.25);
			\draw[line width=1.2pt,orange!88!black,fill=orange!6] (Zc) circle (2);
			\draw[line width=1.2pt,blue!68!black,fill=blue!5] (Zd) circle (1.25);
			\node[font=\small] at (Za) {$8$};
			\node[font=\small] at (Zc) {$5$};
			\node[font=\small] at (Zd) {$8$};
			\node[font=\small,fill=white,inner sep=1.2pt] at (1.15,-3.15) {$-3$};
		\end{scope}
	\end{tikzpicture}
	\caption{The lattice triangle $\mG$ and its oriented Descartes configuration.
		The numbers on the left are squared side lengths; those on the right are
		signed curvatures. The dashed circle, of curvature $-3$, encloses the three
		positive circles.}
	\label{fig:circles}
\end{figure}

Starting from a Descartes configuration and repeatedly replacing one circle
by the other circle tangent to the remaining three produces an Apollonian
packing \cite{LagariasMallowsWilks}. If the initial bends are integers, all subsequent bends are integers, too. The arithmetic and group theory of integral Apollonian packings
were developed systematically by Graham, Lagarias, Mallows, Wilks, and Yan
\cite{GLMWYnumber,GLMWYone}; see also the computational investigations of
Fuchs and Sanden \cite{fuchs} and the exposition \cite{Sarnak}. We recall only the small part needed to connect it to the flea group.

Let $S_i$ be the Descartes reflection replacing the $i$th bend by the other
root of Descartes' equation, and put $\Ap=\langle S_1,S_2,S_3,S_4\rangle$. With
\[
M=\begin{pmatrix}
	1&0&0&0\\0&1&0&0\\0&0&1&0\\1&1&1&-2
\end{pmatrix},
\qquad M(a,b,c,s)^T=(a,b,c,d)^T,
\]
one has $\det M=-2$, so this is a rational rather than a unimodular change
of coordinates. If $P_{ij}$ denotes a coordinate transposition, direct
multiplication gives
\begin{equation}\label{eq:matrix-conjugacy}
	MUM^{-1}=P_{14}S_4,\qquad
	MVM^{-1}=P_{24}S_4,\qquad
	MWM^{-1}=P_{34}S_4.
\end{equation}
Thus a flea jump is a Descartes reflection followed by a relabelling. At the
level of the circle picture, after choosing compatible Euclidean
normalizations, it replaces one circle by the other Soddy circle tangent to
the remaining three.

Let $\Pi_4$ be the group of coordinate permutations and
$\widehat\Ap=\langle\Ap,\Pi_4\rangle$; a superscript $+$ will denote the
determinant-one subgroup. The identities above imply
\begin{equation}\label{eq:group-summary}
	M\Gamma M^{-1}=\widehat\Ap^{+},\qquad
	M\langle U^2,V^2,W^2\rangle M^{-1}=\Ap^+.
\end{equation}
Indeed, the three squares correspond to
$S_1S_4,S_2S_4,S_3S_4$, which generate the even-word subgroup $\Ap^+$,
while the images of $U,V,W$ in the permutation quotient are
$(14),(24),(34)$. Consequently there is an exact sequence
\begin{equation}\label{eq:exact-sequence}
	1\longrightarrow\langle U^2,V^2,W^2\rangle
	\longrightarrow\Gamma\longrightarrow S_4\longrightarrow1.
\end{equation}
The standard Apollonian group is the free product of four groups of order
two, $\Ap\cong C_2*C_2*C_2*C_2$ \cite{GLMWYone}. Hence its even-word
subgroup is a free group of rank $3$. In particular, $\Gamma$ is virtually
free, and the kernel in \eqref{eq:exact-sequence} has index $24$. After
coordinate labels are forgotten, the flea orbit is exactly the
corresponding ordinary Apollonian matrix orbit.

Integral Apollonian packings quickly lead to a local--global question.
Early investigations found strong congruence restrictions on their bends
and suggested that, apart from these restrictions, every sufficiently large
integer should occur \cite{GLMWYnumber,fuchs}. An integer surviving all
congruence tests is called \emph{admissible}. Bourgain and Fuchs proved that
the set of bends has positive density \cite{BourgainFuchs}, and Bourgain and
Kontorovich later proved that almost every admissible integer occurs
\cite{BourgainKontorovich}. These results gave compelling evidence for the
pointwise local--global conjecture \cite{KontorovichSurvey}. The underlying Apollonian group is, however, \emph{thin}: it has infinite index in the integral orthogonal group of the Descartes form while remaining Zariski dense there \cite{KontorovichSurvey,notices}.

In 2024 Haag, Kertzer, Rickards, and Stange showed, however, that congruence conditions do not tell the whole story. They constructed quadratic and quartic reciprocity obstructions and thereby disproved the local--global conjecture for many integral Apollonian packings
\cite{HKRS}. The root $(-3,5,8,8)$ occurring in our flea problem is one of their basic
examples. The quadratic character used is the same
odd-coordinate character from Section~3, and at the root both
computations give $\symb{8}{5}=-1$. However, our no-square theorem is not a direct corollary of their result, because the flea area $|a+b+c-d|/4$ is not itself a bend. This is precisely why Lemma~\ref{lem:ar} is needed.

The contrast is especially visible modulo $24$. In the packing rooted at
$(-3,5,8,8)$, the bends occupy only $0,5,8,12,20,21\pmod{24}$
\cite{HKRS}. Of these classes, only $0$ and $12$ can contain perfect squares,
so the reciprocity obstruction is partly hidden behind ordinary congruence
restrictions. For flea areas, Proposition~\ref{prop:all-residues} says that
every residue class modulo every modulus occurs, yet every square is absent.
Here the global obstruction stands completely exposed.

Had one encountered the flea problem first, the warning would have been hard
to miss. In the actual chronology, the fleas provide a small bridge from
olympiad geometry to a phenomenon discovered only recently in Apollonian
packings.

\noindent\small
\textsc{Giedrius Alkauskas}\\
Vilnius University, Institute of Mathematics and Informatics,\\
Naugarduko 24, LT-03225 Vilnius, Lithuania\\
\href{mailto:giedrius.alkauskas@mif.vu.lt}{giedrius.alkauskas@mif.vu.lt}

\end{document}